\documentstyle[12pt,amscd]{amsart}

\numberwithin{equation}{section}

\newcommand{\CJ}{\cal C(J)}

\newcommand{\hl}{ho\-mot\-o\-py Lie\ }

\newcommand{\nr}{\Bbb{R}}
\newcommand{\nl}{\Bbb{L}}

\newcommand{\bv}{Ba\-ta\-lin-Vil\-ko\-vi\-sky\ }
\def\16N{d\lbrack\phi_1,\dots,\phi_N\rbrack+ \Sigma_{i=1}^N
      \pm \lbrack\phi_1,\dots, d\phi_i,\dots,\phi_N\rbrack=
      \Sigma^{N-1}_{Q=2}\pm \lbrack\lbrack\phi_{i_1},\dots,
      \phi_{i_Q}\rbrack,\phi_{i_{Q+1}},\dots,\phi_{i_N}\rbrack}
\def\L{\cal L}

\def\C2{F(S^1,2)}
\def\C3{\overbar F(S^1,3)}
\def\C4{\overbar F(S^1,4)}

\def\overbar{\overline}

\newcommand{\aaa}{\alpha}
\newcommand{\bb}{\beta}
\newcommand{\ggg}{\gamma}
\newcommand{\dd}{\delta}
\newcommand{\rr}{r^a_\alpha}

\newcommand{\DD}{\Delta}
\newcommand{ \BBvD}{Burgers, Behrends and van Dam}

\newcommand{\smooth}{$C^\infty(M)$}

\newtheorem{th}{Theorem}[section]

\newtheorem{df}{Definition}[section]

\newtheorem{rem}{Remark}[section]

\begin{document}

\title[ascona]{ Deformation Theory and the Batalin-Vilkovisky Master Equation}
\author
{Jim Stasheff $^1$}
\address
{Department of Mathematics, University of North Carolina, Chapel Hill, NC
27599-3250, USA}
\email{jds@@math.unc.edu}
\thanks{$^1$ Research supported in part by NSF grant DMS-9504871. }


\today

\maketitle

\begin{abstract}

The Batalin-Vilkovisky master equations, both classical and quantum,
are precisely the integrability equations for deformations of algebras
and differential algebras respectively.  This is not a coincidence;
the Batalin-Vilkovisky approach is here translated into the language of
deformation theory.
\end{abstract}

The following exposition is based in large part on work by Marc 
Henneaux (Bruxelles) especially and with   Glenn Barnich (Penn State 
and Bruxelles) and Tom Lada and Ron Fulp of NCSU (The Non-Commutative 
State University). The first statement of the relevance of deformation
theory to the construction of interactive Lagrangians, that I am aware
of, is due to Barnich and Henneaux \cite{bh:consis}:
\begin{quote}
We point out that this problem can be economically reformulated as a
deformation problem in the sense of deformation theory \cite{gerst:defm},
namely that of deforming consistently the master equation.
\end{quote}

\vskip2ex

The `ghosts' introduced by  Fade'ev and Popov \cite{fad-pop}
 were soon incorporated into 
the BRST-cohomology approach \cite{BRS} to a variety of problems in 
mathematical 
physics.  There they were reinterpreted by Stora \cite{stora} and others
in terms of the Maurer-Cartan 
forms in the case of a finite dimensional Lie group and more generally as 
generators of the Chevalley-Eilenberg complex \cite{CE} for Lie algebra cohomology.
This led eventually to the Batalin-Vilkovisky approach  
\cite{bv:antired,bv:closure,BV3}
to quantizing particle Lagrangians and then to string field
theory, both classical and quantum \cite{z:csft}.  With hindsight, the
Batalin-Vilkovisky machinery can be recognized as that of homological
algebra \cite {ht}.  The `quantum'  Batalin-Vilkovisky master equation has
the form of the Maurer-Cartan equation for a flat connection, while the
`classical' version has the form of the integrability equation of
deformation theory. In the context of the present conference, my goal is to 
show that these are more than analogies; the master equations are indeed 
the integrability equations of the deformation theory
of, respectively, differential graded commutative algebras and 
graded commutative algebras.

First I will review the jet bundle approach to Lagrangian field theory.
Here we could already encounter cohomological physics in the form of the
variational bicomplex of differential forms on the jet bundle, but I will 
omit that today and proceed instead to the anti-field, anti-bracket 
formalism, the rubric under which physicists reinvented homological algebra.  
Here the `standard construction' is the \bv  complex.

Just as the Maurer-Cartan equation makes sense in the context of Lie
algebra cohomology, so the \bv master equation has an interpretation 
in terms of strong homotopy Lie algebras ($L_\infty$-algebras), as I
explain next. 
After a brief recollection of deformation theory for differential graded
algebras, I will look at various physical examples where free Lagrangians
are deformed to interactive Lagrangians.  Particularly interesting examples are
provided by Zwiebach's closed string field theory \cite{z:csft} 
and higher spin particles \cite{BBvD}.

\section{The jet bundle setting for Lagrangian field theory}

Let us begin with a space  $\Phi$  of {\bf fields} regarded as the space of 
sections of some bundle $\pi:E \to M$.  For expository and coordinate 
computational purposes, I will assume $E$ is a trivial vector bundle and 
will write a typical field as $\phi = (\phi^1,\dots,\phi^k): M\to \nr^k$.
In terms of local coordinates, we start with a trivial vector bundle 
$E= F\times M \to M$ with base manifold M, locally $R^n$, with 
coordinates $x^i, i=1,\dots,n$ and
fibre $R^k$ with coordinates $u^a, a=1,\dots,k$.  We `prolong' this bundle 
to create the associated jet bundle $J=J^\infty E  \to E \to M$ which 
is an infinite dimensional vector bundle with coordinates $u^a_I$ 
where $I=i_1\dots i_r$ is a symmetric multi-index (including, for
$r=0$,  the empty set of indices, meaning just $u^a$).  The notation is 
chosen to bring to mind the mixed partial derivatives of order $r$.  
Indeed, a section of $J$ is the (infinite) jet $j^\infty \phi$
of a section $\phi$ of $E$ if, for all $r$, we have
$\partial_{i_1}\partial_{i_2}...\partial_{i_r} \phi^a=
 u^a_I\circ j^\infty \phi$ where $\phi^a = u^a\circ \phi$ and $\partial_i =
\partial/\partial x^i$.
\begin{df}
A {\bf local function} $L(x,u^{(p)})$ is a smooth function in the
coordinates $x^i$ and the coordinates $u^a_I$, where the order $|I| = r$ of the
multi-index $I$ is less than or equal to some integer $p$.
\end{df}

Thus a local function is in fact the pullback of a smooth function on some
finite jet bundle $J^pE$, i.e.  a composite $J \to J^pE \to R$.

The {\bf space of local functions} will be denoted $\CJ$.
\begin{df}
A {\bf local functional}
\begin{eqnarray}
{\cal L}[\phi]=\int_M  L(x,\phi^{(p)}(x)) d{vol}_M
= \int_M (j^\infty \phi)^*  L(x,u^{(p)}) d{vol}_M
\end{eqnarray}
is the integral over $M$ of a local function evaluated for 
sections $\phi$ of $E$.
(Of course, we must restrict $M$ and $\phi$ or both for this to make sense.)
\end{df}

The variational approach is to seek the critical points of such
a local functional. More precisely, we seek sections $\phi$ such
that $\delta {\cal L}[\phi]=0$ where $\delta$ denotes the
variational derivative corresponding to an `infinitesimal'
variation: $\phi \mapsto \phi +\delta\phi$.  The condition
 $\delta {\cal L}[\phi]=0$ is equivalent to the Euler-Lagrange
equations on the corresponding local function $L$ as follows:
Let
\begin{eqnarray}
D_i={\partial\over\partial x^i}+u^a_{Ii}{\partial\over\partial
u^a_{I}}
\end{eqnarray}
be the total derivative acting on local functions and
\begin{eqnarray}
E_a=(-D)_I{\partial\over\partial u^a_{I}}
\end{eqnarray}
the Euler-Lagrange derivatives.  The notation $(-D)_I$ means
$(-1)^r D_{i_1}\cdots D_{i_r}$.
The Euler-Lagrange equations are then
\begin{eqnarray}
E_a(L)=0.
\end{eqnarray}
 
Two local functionals ${\cal L}$ and ${\cal K}$ are equivalent
if and only if the
Euler-Lagrange derivatives of their integrands agree, $E_a(L)=E_a(K)$,
for all $a=1,\dots,k$.
 The kernel of the Euler-Lagrange derivatives is given by total
divergences,
\begin{eqnarray}
E_a(L)=0,\ for\ all\ a=1,\dots,k \Longleftrightarrow L=(-1)^k D_i j^i
\end{eqnarray}
for some local functions $j^i$. Equivalently,
\begin{eqnarray}
E_a(L)=0,\ for\ all\ a \Longleftrightarrow L(x,u^{(p)})
d{vol}_M = (-1)^{i-1}d j^i  d{vol}_M/d x^i 
\end{eqnarray}
where if
 $d{vol}_M = dx^1\cdots dx^n,$ 
then $d{vol}_M/d x^i = dx^1\cdots dx^{i-1} dx^{i+1} \cdots dx^n$.

Since $\cal L$ is the integral of an $n$-form on $J$, it is not 
surprising that this all makes sense in the deRham complex
$\Omega^*(J)$,
which remarkably splits as a bicomplex (though the finite level
complexes $\Omega^*(J^pE)$ do not).  The appropriate 1-forms in the
fibre directions are not the $du^a_I$ but rather the {\bf contact forms}
$\theta^a_I= du^a_I- u^a_{Ii}dx^i $.  The total differential $d$ splits as
$d=d_H+d_V$ where $d_H = dx^i\partial/\partial x^i$ is the usual
exterior differential on $M$ pulled up to $J$.  We will henceforth
restrict the coefficients of our forms to be local functions, although we
will not decorate $\Omega^*(J)$ to show this.

The Euler-Lagrange operators assemble into an operator on forms:
\begin{eqnarray}
E(L dvol_M) = E_a(L)\theta^a dvol_M.
\end{eqnarray}
A Lagrangian $\cal L$  determines a {\bf stationary surface} 
or   {\bf solution surface} or 
{\bf shell} $\Sigma\subset J^\infty$ such that $\phi$ 
is a solution of the variational problem (equivalently,
the Euler-Lagrange equations)
 if $j^\infty \phi$ has its image in $\Sigma$.  
The corresponding algebra is the {\bf stationary ideal} $\cal I$
of local functions which vanish `on shell', i.e. when resticted to
the  solution surface $\Sigma$.

The Euler-Lagrange equations generate $\cal I$ as a differential ideal,
but this means we may have not only {\bf Noether identities}
\begin{eqnarray}
r^a_\alpha E_a(L) = 0
\end{eqnarray}
but also 
\begin{eqnarray}
r^{a I}_\alpha  D_IE_a(L) = 0.
\end{eqnarray}
Of course we have `trivial' identities of the form 
\begin{eqnarray}
D_JE_b(L)\mu_\alpha^{bJaI}D_IE_a(L)=0,
\end{eqnarray}
since we are dealing with a commutative algebra of functions.
One can show
\cite{ht} that all Noether identities, where the coefficients $r^a_\alpha$
vanish on shell, are of the above form.
We now assume we have a set of indices $\{\alpha\}$ such that the above
identities generate all the non-trivial relations in $\cal I$.  
According to Noether
\cite {noether}, each such identity corresponds to 
an {\bf infinitesimal  gauge symmetry}, i.e. an infinitesimal variation
that preserves the space of solutions or, equivalently, a vector field
tangent to $\Sigma$. For each Noether identity indexed by $\alpha$, 
we denote the
corresponding vector field by $\delta_\alpha$.
We denote by 
 $\Xi$, the 
{\bf space  of gauge symmetries}, considered as a vector space but
also  as a module over $\CJ$.
We can regard $\delta_\alpha$ as a  
(constant) vector field on the space of fields $\Phi$ and hence $\delta$ 
as a linear map
$$
\delta:\Xi\to Vect\Phi.
$$

Since the bracket of two such vector fields $\lbrack\delta_\alpha,\delta_
\beta\rbrack$ is again a gauge symmetry, it agrees with something in the 
image of $\delta$ when acting on solutions.  If we denote that something
as $\lbrack\alpha,\beta\rbrack$, one says this bracket `closes on shell'.  
It is not in general a Lie bracket, since 
the Jacobi identity may hold only `on shell'.  
(Later I will address the issue that 
$\lbrack\delta_\alpha,\delta_ \beta\rbrack$ 
may not be constant on $\Phi.$) 

To make this more explicit, write
\begin{eqnarray}
\lbrack\delta_\alpha,\delta_\beta\rbrack=
\delta_{\lbrack\alpha,\beta\rbrack} + \nu^{a}_{\alpha\beta}{{\delta L} 
\over {\delta u^a}}.
\end{eqnarray}
The possible failure of the Jacobi identity results from those last terms 
which vanish only on shell and the fact that we are working in a module 
over $\CJ$. (For example, we have structure {\it functions} rather 
than structure constants in terms of our generators.)

All of this, including these latter subtleties, are incorporated into
a remarkable complex by \bv, which we shall describe below.

\section{Deformation theory}

A deformation theoretic approach to producing interactive Lagrangians is 
to start with a `free' Lagrangian $L_0$ on a space of fields $\Phi$ with an 
abelian Lie algebra  $\Xi$ represented as gauge symmetries by a Lie 
map 
$$
\delta_0:\Xi\to Vect\Phi,
$$
and then study formal deformations
$$L = \Sigma t^i L_i,\ \  \delta = \Sigma t^i \delta_i,$$
keeping $\Phi$ and $\Xi$ fixed as vector spaces,
such that $\delta_\xi\int L dvol_M = 0$ for all $\xi\in\Xi.$

Classical deformation theory emphasizes the following two problems:

$\bullet$ 1.  Consider a candidate {\bf infinitesimal} $L_1$ and see if there 
exists a full formal deformation $L$ with corresponding $\delta$.  Such candidate 
terms $L_1$ are often suggested by physical descriptions of elementary 
interactions.

$\bullet$ 2. Classify all formal deformations up to appropriate equivalence.

A third subsequent problem, that of convergence of the power series 
involved, is a problem in analysis; cohomological techniques apply to the 
formal algebraic theory - I don't do estimates.

The cohomological approach to deformation theory, as initiated by 
Gerstenhaber \cite{gerst:defm}, situates the problem in an appropriate complex, the 
Hochschild cochain complex $C^*(A,A)$ in the case of an associative 
algebra $A$.  For the Lagrangian problem, the complex is due to Batalin and 
Vilkovisky \cite{bv:antired,bv:closure,BV3}
 using anti-field and ghost technology and the anti-bracket 
of Zinn-Justin \cite{zj}.

\section{Anti-fields, (anti-)ghosts and the anti-bracket}
	Let me take you `through the looking glass' and present a `bi-lingual' 
(math and physics) dictionary.

{} From here on, we will talk in terms of algebra extensions 
of $C^\infty (E)$ and $\CJ$, 
but the extensions will all be free graded commutative. 
We could instead talk in terms of an extension of $E$ or $J$ 
as a super-manifold, 
the new generators being thought of as (super)-coordinates.  

We first  extend $C^\infty (E)$ by adjoining generators of 
various degrees to form a free graded commutative algebra $\cal A_1$ 
over $C^\infty (E)$, that is, even graded generators 
give rise to a polynomial algebra and odd graded generators give rise to a 
Grassmann ($=$ exterior) algebra.  
The generators (and their 
products) are, in fact, bigraded $(p,q)$; the graded commutativity
is  with respect to the total degree $p-q$.

For each variable $u^a$, adjoin an anti-field $u_a^*$ 
and for each $r_\alpha$,  adjoin a corresponding ghost $C^\alpha$ 
and a corresponding anti-ghost $C^*_\alpha$. Here is a table showing the 
corresponding math terms and the bidegrees.
\vskip4ex
\begin{tabular}{|l|l|c|c|c|}\hline
Physics & Math & Ghost & Anti-ghost & Total\\
Term & Term & Degree & Degree & Degree \\
\hline
field & section & 0 & 0 & 0 \\
\ & \ & \ & \ & \ \\
anti-field & Koszul generator & 0 & 1 & -1 \\
\ & \ & \ & \ & \ \\
ghost & Cartan-Eilenberg  generator & 1 & 0 & 1 \\
\ & \ & \ & \ & \ \\
anti-ghost & Tate  generator & 0 & 2 & -2 \\
\hline
\end{tabular}

\vskip4ex
Note that the anti-field coordinates depend on $E$ but the ghosts and 
anti-ghosts depend also on the specific Lagrangian.

This algebra in turn can be given an {\bf  anti-bracket} $(\ ,\ )$ of degree 
$-1$ which, remarkably, combines with the product we began with to 
produce precisely an `up-to-homotopy' analog of a
  Gerstenhaber algebra \cite{lz:new,kvz}, 
though this was not 
recognized until quite recently.  

\begin{df} A {\bf Gerstenhaber algebra} is a graded commutative and
associative algebra $A$ together with a bracket $[\cdot,\cdot]: A\otimes
A\,\to\, A$
of degree $-1$, such that for all homogeneous elements $x$, $y$, and $z$ in
$A$,
\[
[x,y] := -(-1)^{(|x|-1)(|y|-1)} [y,x],
\]
\[
[x,[y,z]] = [[x,y],z] + (-1)^{(|x|-1)(|y|-1)} [y,[x,z]],
\]
and
\[
[x, y z] = [x,y] z + (-1)^{(|x|-1)|y|} y [x,z].
\]
\end{df}

In the field-anti-field formalism (without jet coordinates),
the anti-bracket looks like
\begin{eqnarray}
(\phi^a(x), \phi^*_b(y)) = \delta^a_b \delta (x-y)
\end{eqnarray}
where $\delta (x-y)$ is the Dirac delta `function' (distribution).
In the corresponding jet bundle formalism,
the anti-bracket is defined on generators 
and then extended to polynomials by applying the graded Leibniz identity 
so that $(\psi,\ )$ is a graded derivation for any $\psi$ in this algebra.
The only non-zero anti-brackets of generators are
$$
(u^a,u_b^*) = \delta^a_b \ \ {\text {and}\ \ } (C^\alpha, C^*_\beta) 
=\delta^\alpha_\beta.$$

Now we further extend  $\CJ$ with corresponding jet 
coordinates $u^{I*}_a, C^\alpha_I$ and $C^{I*}_\alpha$ with the 
corresponding pairings giving the extended anti-bracket. The
resulting \bv algebra we denote $\cal A$.

\def\L{L_\infty}
\def\LL{$$\overbar L = L + \phi^*_ir^i_\alpha C^\alpha 
+ (C_{\ggg}^*c^{\ggg}_{\alpha\bb} + \phi^*_i\phi^*_j\4mu4) C^\alpha C^\bb$$}
\def\SS{$\overbar S =\int_M\overbar L dvol_M$}
\def\s{s=(\overbar L,\ \ ) }
\def\umu{u^{i\overbar \mu}}
\def\4mu4{\nu^{a}_{\alpha\beta}}

Define an operator $s_0$ of degree $-1$ on $\cal A$ as $ (L,\ )$.

We call the antifields Koszul generaters because
\begin{eqnarray}
s_0 u^*_a =  {{\dd L}\over {\dd u^a}} 
\end{eqnarray}
as in the Koszul complex for the ideal,
so that $H^{0,0}\subset \Phi$ is given by  
${{\dd L}\over {\dd u^a}} = 0,$
but   
\begin{eqnarray}
s_0(r^a_\alpha u^*_a) = r^a_\alpha  {{\dd L}\over {\dd u^a}} =0,
\end{eqnarray}
so that $H^{0,1} \neq 0.$

Now consider the extended Lagrangian 
\begin{eqnarray}
L_1 = L_0 + u^*_a r^a_\alpha C^\alpha
\end{eqnarray}
and $s_1 = (L_1,\ ),$
so that 
\begin{eqnarray}
s_1 C^*_\alpha = u_a r^a_\alpha
\end{eqnarray}
and $H^{0,1}$ is now $0$, as in Tate's extension of the Koszul
complex of the ideal to produce a resolution \cite{tate}.  That is why
we refer to the anti-ghosts as Tate generators.
(If needed, Tate tells 
us to add further generators in bidegree $(0,q)$ for $q >2$ so that $H^{0,q} 
= 0$ for $q>0.$)

\def\cc{ c^\alpha_{\bb\ggg}C^\bb C^{\ggg}}

Further extend $L_1$ to
\begin{eqnarray}
L_2 = L_1 + C^\aaa_*\cc,
\end{eqnarray}
so that
$$
s_2C^\alpha = c^\alpha_{\bb\ggg}C^\bb C^{\ggg}
$$
$$
s_2u^a = \rr C^\alpha
$$
which is how the Chevalley-Eilenberg coboundary looks in terms of bases
for a Lie algebra and a module and corresponding structure constants.  
However, we may not have $(s_2)^2 = 0$ since $\rr$ and $c^\alpha_{\bb\ggg} $ 
are functions.
Batalin and Vilkovisky prove that all is not lost.  First, they add to
$L_2$ a term involving the functions $\4mu4$.

\begin{th} $L_2$ can be further extended by terms of higher degree in the 
anti-ghosts to $L_\infty$  so that $(L_\infty,L_\infty) = 0$ and
hence the corresponding $s_\infty$ will have square zero.
\end{th}

With hindsight, we can see that the existence of these terms of higher 
order is guaranteed because the antifields and antighosts provide a 
resolution of the stationary ideal.

We refer to this complex $(\cal A, s_\infty)$ as the {\bf \bv complex}.

What is the significance of $(s_\infty)^2 = 0$ in our Lagrangian
context,  or, equivalently, of the {\bf Master Equation} $(\L,\ L) = 0$?  
There are three answers:  in higher 
homotopy algebra, in deformation theory and in mathematical physics.
It is the deformation theory that provides the transition betweeen
the other two.

\section{The Master Equation and Higher Homotopy Algebra}

If we expand $s = s_0 + s_1 + ...$ where the subscript indicates the change 
in the ghost degree $p$, the individual $s_i$ do not correspond 
to $(x, \ )$ for 
any term $x$ in $\L$ but do have the following graphical description:
\vskip2ex
standard diagram of differentials of a spectral sequence goes here
\vskip2ex
\noindent so that we see the B-V-complex as a multi-complex. 
The differential $s_1$
gives  us  the Koszul-Tate differential $d_{KT}$
and part of $s_2$ looks like that of 
Chevalley-Eilenberg.  That is, $C_\alpha^* c^\alpha_{\bb\ggg}C^\bb C^{\ggg}$ describes the
(not-quite-Lie) bracket on $\Xi$.  Further terms with one anti-ghost 
$C_\alpha^*$ and three ghosts $C^\bb C^{\ggg} C^\delta$ describe a tri-linear $[\ ,\ ,\ ]$ 
and so on for multi-brackets of possibly arbitrary length.  Moreover, the 
graded commutativity of the underlying algebra of the B-V-complex implies
appropriate symmetry of these multi-brackets.  The condition that 
$s_\infty^2=0$ translates to the following identities, which are the defining 
identities for a {\bf strong homotopy Lie algebra} or  $L_\infty$ 
{\bf algebra}.  
\begin{eqnarray}
d[v_1, \dots, v_n] + \sum_{i=1}^n  \epsilon(i) [v_1, \dots, dv_i, \dots, v_n]
\nonumber\\
= \sum_{ k+l = n+1}
\sum_{ \text{unshuffles } \sigma}
\epsilon (\sigma) [[v_{i_1}, \dots, v_{i_k}], v_{j_1}, \dots, v_{j_{l}}],
\end{eqnarray}
where $\epsilon (i) = (-1)^{|v_1| + \dots + |v_{i-1}|}$ is the sign picked up
by taking $d$ through $v_1,  \dots,  v_{i-1}$ and, for the unshuffle
$\sigma:\{1,2, \dots, n\}\mapsto  \{i_1, \dots, i_k, j_1, \dots, j_{l}\},$ 
the sign $\epsilon (\sigma)$ is the sign
picked up by the elements $v_i$ passing through the $v_j$'s during the
unshuffle of $v_1, \dots , v_n$, as usual in superalgebra.

\begin{rem}
Here we follow the physics grading and  sign conventions in our definition of
a strong \hl algebra \cite{wz,z:csft}. These are equivalent to but different 
from those
in the existing mathematics literature, cf.\ Lada and Stasheff \cite{ls},
in which the $n$-ary bracket has degree $2-n$. (The correspondence
and additional insights are presented in full detail in Kjeseth's
dissertation \cite{lars}.)  With those mathematical
conventions, $L_\infty-$algebras occur naturally as 
deformations of Lie algebras.
If $L$ is a Lie algebra and $V$ is a complex with a homotopy equivalence to the
trivial complex $0 \to L \to 0$, then $V$ is naturally a \hl algebra, see
Schlessinger and Stasheff \cite{SS}, Barnich, Fulp, Lada and Stasheff
\cite{bfls} and  Getzler and Jones \cite{gj:n-bv}.
\end{rem}

Realized in the \bv complex,
 these defining identities tell us, for small values of $n$, that $d_{KT}$ 
is a graded derivation of the bracket, that the bracket may not satisfy the 
graded Jacobi identity but that we do have (with the appropriate signs)
\begin{multline}
[[v_1, v_2], v_3] \pm [[v_1, v_3], v_2] \pm [[v_2, v_3], v_2] =
\\
- d[v_1, v_2, v_3] \pm [dv_1, v_2, v_3]\pm[v_1, dv_2, v_3]\pm [v_1, v_2,
dv_3].
\end{multline}
i.e. the Jacobi identity holds {\it up to homotopy} or, for closed forms, the 
Jacobi identity holds modulo an exact term - the tri-linear bracket.

Note that the identity is the Jacobi identity if $d_{KT}=0$ and all the 
other brackets vanish and that the identity has content even if only one 
$n$-linear bracket is non-zero and all the others vanish.  Precisely that 
situation has recent been studied quite independently of my work and of 
each other by Hanlon and Wachs \cite{hanlon-wachs} (combinatorial algebraists), 
by Gnedbaye \cite{gned}
(of Loday's school) and by Azcarraga and Bueno \cite{az-bu} (physicists).  
At the Ascona conference, Flato brought to my attention that Takhtajan's 
identity \cite{takh:nambu} for his trilinear `bracket' is not the one above, 
but suitably symmetrized does agree with it. (This identity was also
known to Flato and Fronsdal in 1992, though unpublished.)

\section{The Master Equation and Deformation Theory}

The Master Equation $(L_\infty, L_\infty) = 0$ has precisely the form
of Gerstenhaber's condition for $\L$ to be a deformation of $L_0$. 
Classical (formal) algebraic deformation theory uses a differential graded
Lie algebra (dgla) $\nl$ (e.g. the Hochschild cochain complex with
the Gerstenhaber bracket) to study the problem of `integrating' an
infinitesimal deformation $\theta$ to a full formal deformation $\theta_t$.
The primary obstruction, regarding $\theta$ as a class in $H^*$ of this
dgla $\nl$,
is $[\theta,\theta]$ and further obstructions can be described in terms of
multi-brackets on this cohomology.  Alternatively, the formal deformation
$\theta_t$ itself as a cochain must satisfy $[\theta_t,\theta_t]=0$.
The analogy with the Master Equation is manifest.

Following the historical pattern in algebraic deformation theory,
we could hope to calculate this homology  to be
$0$ in the relevant dimensions  in certain cases, thus obtaining results
of unobstructedness for the integrability question or of rigidity for the
classification problem.  Such
calculations are highly non-trivial, however, and to my knowledge have
been carried out only in the case of  electricity and magnetism
(Maxwell's equations), Yang-Mills and gravity.

\section{The Master Equation in Field Theory}

In the Lagrangian setting, we wish to deform not just the  local
functional, but rather the underlying local function $L$. 
In the case of  electricity and magnetism
(Maxwell's equations), Yang-Mills and gravity, 
 the relevant algebra of gauge symmetries is described by a finite 
dimensional Lie algebra which, moreover, holds off shell.
In terms of an appropriate basis and in the notation of section 1, we
have
\begin{eqnarray}
\lbrack\delta_\alpha,\delta_\beta\rbrack=
c^\ggg_{\aaa\bb}\delta_{\gamma}
\end{eqnarray}
for structure {\it constants} $c^\gamma_{\alpha\beta}$ and acting
on all fields, not just on solutions. This allows the extended Lagrangian 
to be no more than quadratic in the ghosts.  

As field theories, the electron can be described by a field of spin 1,
as can a Yang-Mills particle, while  the graviton   can be described by
a field of spin 2.  Somehow this is related to the strict Lie algebra
structures just described. For higher spin particles, however, we have 
quite a different story, 
which first caught my attention in the work of Burgers, Behrends and van 
Dam \cite{BBvD,burgers}, though I have since learned there was quite a 
history before that and 
major questions still remain open.  By higher spin particle Lagrangians, I 
mean that 
the fields are symmetric $s$-tensors (sections of the symmetric $s$-fold 
tensor product of the tangent bundle).  
If the power is $s$, the field
is said to be of spin $s$ and represents a particle of spin $s$. 
\BBvD \ start with a free theory with 
abelian gauge symmetries and calculate all possible infinitesimal 
interaction terms up to the appropriate equivalence (effectively 
calculating the appropriate homology group).  They then sketch the 
problem of finding higher order terms for the Lagrangian, but do not carry 
out the full calculation.  In fact, according to the folklore in the subject,
a consistent theory for $s\geq 3$ will require additional fields of 
arbitrarily high spin $s$.  For $s=3$, the conjecture is that all higher 
integral spins are needed.  From the deformation theory point of view, this 
suggests the following attack:  Compute the primary obstructions and 
discover that all infinitesimals are obstructed.  Add additional fields to 
kill the obstructions and calculate that indeed additional fields of 
arbitrarily high spin $s$ are needed.  In one memorable phrase, this would 
be `doing string field theory the hard way''.

Zwiebach \cite{z:csft} does indeed have a consistent closed string field 
theory (CSFT),
but produced in an entirely different way.  Recall one of the earliest 
examples of deformation quantization, the Moyal bracket.  Moyal was able 
to produce a non-trivial deformation of a commutative algebra \smooth  
\ on a symplectic manifold, with infinitesimal given by the  Poisson bracket,
by writing down the entire formal power series.  Similarly, Zwiebach is 
able to describe the entire CSFT Lagrangian (at tree level) by giving it in 
terms of the differential geometry of the moduli space of punctured 
Riemann spheres (tree level $=$ genus $0$).  In fact, Zwiebach has the 
following structure:
a differential graded Hilbert space $(\cal{ H}, <\ ,\ >, Q)$ related to the 
geometry of 
the moduli spaces from which he deduces $n$-ary operations $[\ ,\dots\ ,\ ]$ 
which give an $L_\infty$ structure:
$$
\16N
$$
\noindent with $[\phi] = Q\phi$.  

The deformed Lagrangian (still classical) and hence the Master Equation
is satisfied  for
$$S(\Psi ) = {1\over 2} \langle \Psi , Q\Psi \rangle
+ \sum_{n=3}^\infty {\kappa^{n-2} \over n!}
\{ \Psi\dots\Psi \}.
$$
The expression $\{ \Psi\dots\Psi \}$ contains n-terms and will be abbreviated
$\{ \Psi^n \};$ it is given in terms of the brackets by
$\{ \Psi\dots\Psi \} = \langle \Psi , [ \Psi,\dots, \Psi ] \rangle.$
\par
The field equations follow from the classical action by simple variation:
$$\delta S = \sum_{n=2}^\infty {\kappa^{n-2} \over n!} 
\{ \delta\Psi , \Psi^{n-1} \}
$$

with gauge symmetries given by

$$\delta_\Lambda \Psi 
= \sum_{n=0}^\infty { \kappa^n \over n!} [\Psi^n, \Lambda ] .$$
 
\section{Quantization}

So far our description of the anti-field, anti-bracket formalism has been 
in the context of deformations of `classical' Lagrangians.  Batalin and 
Vilkovisky (as well as much of the work on BRST cohomology) were 
motivated by problems in quantum therory.  The quantum version of the 
anti-field, anti-bracket formalism involves a further `second order' 
differential operator $\DD$ of square $0$ on the B-V complex relating the 
graded commutative product and the bracket - namely, the bracket is the 
deviation of the operator $\DD$ from being a derivation of the product.  
This has led to the abstract definition of a BV-algebra.

\begin{df} A {\bf BV-algebra} is a Gerstenhaber algebra with an 
operator (necessarily 
of degree $-1$ for a bracket of degree $-1$) such that 
$$
[A,B] = \DD(AB) -\DD(A)B +(-1)^A\DD(B).
$$
\end{df}

Alternatively, a definition can be given in terms of a graded commutative 
algebra with an appropriate operator $\DD$ \cite{akman:diffop,schw:bvgeom}.

The quantization of Zwiebach's CSFT involves further expansion of the 
Lagrangian in terms of (the moduli space of) Riemann surfaces of genus 
$g\geq 0.$  Here the operator $\DD$ is determined by the self-sewing of a 
pair of pants (a Riemann sphere with 3 punctures).  Now Zwiebach's CSFT
provides a solution of the `quantum Master Equation' which
in the context of a BV-algebra, is
$$
(S,S)=\DD S.
$$
Again we see an anolog of the Maurer-Cartan equation or of a flat 
connection, but why?
\vskip3ex


\end{document}